# ISS Estimates in the Spatial Sup-Norm for Nonlinear 1-D Parabolic PDEs

**Iasson Karafyllis[*] and Miroslav Krstic[**]**

[*]Dept. of Mathematics, National Technical University of Athens,
Zografou Campus, 15780, Athens, Greece, email: iasonkar@central.ntua.gr

[**]Dept. of Mechanical and Aerospace Eng., University of California, San Diego, La Jolla, CA 92093-0411, U.S.A., email: krstic@ucsd.edu

**Abstract**
This paper provides novel Input-to-State Stability (ISS)-style maximum principle estimates for classical solutions of highly nonlinear 1-D parabolic Partial Differential Equations (PDEs). The derivation of the ISS-style maximum principle estimates is performed by using an ISS Lyapunov Functional for the sup norm. The estimates provide fading memory ISS estimates in the sup norm of the state with respect to distributed and boundary inputs. The obtained results can handle parabolic PDEs with nonlinear and non-local in-domain terms/boundary conditions. Three illustrative examples show the efficiency of the proposed methodology for the derivation of ISS estimates in the sup norm of the state.

**Keywords:** ISS, parabolic PDEs, boundary disturbances.

## 1. Introduction

Rapid progress in the development of Input-to-State Stability (ISS) theory for systems modeled by Partial Differential Equations (PDEs) took place during the last decade (see for instance [2,6,7,8,11,13,15,19,20,21,25,26]). Researchers dealt with the two major problems that arise in the study of PDEs with inputs and do not appear in the study of finite-dimensional systems or delay systems: (i) the selection of the state norm (functional norms are not equivalent), and (ii) the presence of boundary inputs (the boundary inputs enter through unbounded operators). Many methodologies have been used for the derivation of ISS estimates (Lyapunov functionals, spectral methods, semigroup methods, etc.). The development of ISS theory for PDEs has produced novel results which can be used for various purposes in control theory and practice. ISS theory for PDEs has been used for the stability analysis of coupled PDEs (see [2,15]), for feedback design purposes in linear and semilinear PDEs (see [12,14,24]) and for observer design purposes (see [16]) for various infinite-dimensional systems.

   Although the sup norm is the most useful norm (because it can provide point-wise estimates), there are few ISS estimates in the literature for the sup norm of the solution of parabolic PDEs. The scarcity of ISS estimates in the state sup norm is justified by the absence of ISS Lyapunov functionals which are related to the sup norm. ISS estimates in the sup norm were provided in the work [11], where numerical approximations of the solution of a linear 1-D parabolic PDE were exploited. The methodology was later generalized in [15] by introducing the notion of an ISS Lyapunov functional under discretization. Input-to-Output Stability (IOS) estimates for the state sup norm were provided for certain semilinear parabolic PDEs in [26]. Indeed, in [26] the state sup norm was estimated by using the $H^1(0,1)$ norm of the initial condition and consequently, the provided estimates were IOS estimates with output map being the state in the space $L^\infty(0,1)$ while the state space was the Sobolev space $H^1(0,1)$. ISS estimates of the $H^1(0,1)$ norm of the state for certain linear 1-D parabolic PDEs were also provided in [13]. These estimates can be used for the derivation of estimates of the sup norm by using auxiliary tools (e.g., Agmon's inequality).

The present work provides ISS estimates in the sup norm for highly nonlinear and non-local 1-D parabolic PDEs by using a completely different methodology than the methodology used in [11,15]. Instead of exploiting numerical approximations of the solution of the PDE problem, here we construct an ISS Lyapunov Functional for the sup norm. To that purpose, we consider classical solutions of the 1-D parabolic PDE (as in [11,15] but less regular solutions than the solution notion in [11,15]). The constructed ISS Lyapunov functional is a coercive Lyapunov functional (i.e., is bounded from above and below by $K_\infty$ functions of the sup norm of the state). It should be noticed that non-coercive Lyapunov functionals may be used for the proof of the ISS property (see [8]). Using the coercive ISS Lyapunov functional, we are able to get more general results than the results in [11,15], which can deal with highly nonlinear (and non-local) parabolic PDEs and highly nonlinear (and non-local) boundary conditions. More specifically, we obtain ISS-style maximum principle estimates which (like standard maximum principles; see [3,9,22]) provide bounds for the sup norm of the state that depend only on the sup norm of the initial condition, the source terms in the PDE and the boundary values of the solution (and its spatial derivative). However, unlike standard maximum principles, the obtained estimates take into account the exponential decay of the effect of the initial condition and the fact that recent input values have stronger effect on the current value of the solution than past input values (fading memory effect). The ISS-style maximum principle estimates provide in a direct way fading memory ISS estimates when disturbance inputs appear in the boundary conditions. The proof methodology is simpler and less lengthy than the methodology used in [11,15].

The structure of the paper is as follows. Section 2 contains the statements of the main results of the paper and provides certain remarks. The main results are used in Section 3 in three (carefully chosen) illustrative examples, which show how easily the main results can be applied to highly nonlinear parabolic PDEs. The proofs of the main results are provided in Section 4, where some interesting auxiliary results are also stated and proved. Section 5 gives the concluding remarks of the present work.

**Notation.** Throughout this paper, we adopt the following notation.

∗ Let $U \subseteq \Re^n$ be a set with non-empty interior and let $\Omega \subseteq \Re$ be a set. By $C^0(U;\Omega)$, we denote the class of continuous mappings on $U$, which take values in $\Omega$. By $C^k(U;\Omega)$, where $k \geq 1$, we denote the class of continuous functions on $U$, which have continuous derivatives of order $k$ on $U$ and take values in $\Omega$. When $\Omega \subseteq \Re$ is not explicitly given, i.e., when we write $C^k(U)$, we mean that $\Omega = \Re$. For $\eta \in C^1([0,1])$, $\eta'(x)$ denotes the derivative with respect to $x \in [0,1]$.

∗ $\Re_+$ denotes the set of non-negative real numbers. A function $\rho : \Re_+ \to \Re_+$ is said to be positive definite if $\rho(s) > 0$ for all $s > 0$ and $\rho(0) = 0$. A positive definite function $\rho \in C^0(\Re_+;\Re_+)$ is said to be a function of class $K_\infty$ if $\rho$ is increasing with $\rho(\Re_+) = \Re_+$. The sign function $\text{sgn}(x)$ is defined by $\text{sgn}(x) = x/|x|$ for $x \neq 0$ and $\text{sgn}(0) = 0$.

∗ Let $u \in C^0(I \times [0,1])$ be given, where $I \subseteq \Re$ is an interval. We use the notation $u[t]$ to denote the profile at certain $t \in I$, i.e., $(u[t])(x) = u(t,x)$ for all $x \in [0,1]$. $L^\infty(0,1)$ denotes the equivalence class of Lebesgue measurable functions $f : [0,1] \to \Re$ for which $\|f\|_\infty = \underset{x \in (0,1)}{\text{ess sup}}(|f(x)|) < +\infty$.

## 2. ISS in the Spatial Sup-Norm

Our first main result provides an estimate of a weighted sup norm of the solution of a PDE with space and time-varying coefficients. The estimate depends only on the initial condition, the distributed perturbation term that may be present in the PDE and the boundary values of the solution.

**Theorem 2.1:** *Let $T > 0$, $a, b, c : (0, T) \times (0,1) \to \Re$, $f : (0, T] \times (0,1) \to \Re$ with $\sup_{0 < x < 1, 0 < t \leq T} (|f(t, x)|) < +\infty$, $u \in C^0([0, T] \times [0,1])$ with $u[t] \in C^2((0,1))$ for almost all $t \in (0, T)$, for which the derivative $\frac{\partial u}{\partial t}(t, x)$ exists, is a continuous function on $(0, T) \times [0,1]$ and satisfies*

$$\frac{\partial u}{\partial t}(t, x) = a(t, x) \frac{\partial^2 u}{\partial x^2}(t, x) + b(t, x) \frac{\partial u}{\partial x}(t, x) + c(t, x) u(t, x) + f(t, x),$$

*for almost all $t \in (0, T)$ and for all $x \in (0,1)$*  (2.1)

*Suppose that*

$$a(t, x) \geq 0, \text{ for all } t \in (0, T), \ x \in (0,1) \tag{2.2}$$

*and that there exist a constant $\sigma > 0$ and a positive function $\eta \in C^2([0,1]; (0, +\infty))$ such that*

$$a(t, x) \eta''(x) + b(t, x) \eta'(x) + (\sigma + c(t, x)) \eta(x) \leq 0, \text{ for all } t \in (0, T), \ x \in (0,1) \tag{2.3}$$

*Then the following estimate holds for all $\zeta \in [0, \sigma)$ and $t \in (0, T]$:*

$$\|u[t]\|_{\infty, \eta} \leq \max \left( \exp(-\zeta t) \|u[0]\|_{\infty, \eta}, \sup_{0 < s \leq t} \left( \max \left( \frac{|u(s, 0)|}{\eta(0)}, \frac{|u(s,1)|}{\eta(1)}, \frac{\|f[s]\|_{\infty, \eta}}{\sigma - \zeta} \right) \exp(-\zeta(t-s)) \right) \right)$$
(2.4)

*where*

$$\|u[t]\|_{\infty, \eta} := \max_{0 \leq x \leq 1} \left( \frac{|u(t, x)|}{\eta(x)} \right), \ \|f[t]\|_{\infty, \eta} := \sup_{0 < x < 1} \left( \frac{|f(t, x)|}{\eta(x)} \right). \tag{2.5}$$

Estimate (2.4) is an ISS-style maximum principle estimate which (like ordinary maximum principles; see [3,9,22]) provides a bound for the sup norm of the state that depends only on the sup norm of the initial condition, the source terms in the PDE ($f$) and the boundary values of the solution. However, unlike ordinary maximum principles, estimate (2.4) takes into account the exponential decay of the effect of the initial condition and the fact that past input values have less effect on the current value of the solution than more recent input values (fading memory effect). The ISS-style maximum principle estimate (2.4) provides in a direct way a fading memory ISS estimate when Dirichlet boundary inputs are present, i.e. when we have the boundary conditions

$$u(t, 0) = d_0(t), \ u(t, 1) = d_1(t), \text{ for all } t \in [0, T]$$

where $d_0, d_1 \in C^0([0, T])$ are the boundary disturbance inputs.

The proof of Theorem 2.1 relies on showing that the weighted sup norm of the solution $V(u) = \|u\|_{\infty, \eta}$ (defined by (2.5)) is an ISS Lyapunov functional for the PDE system (2.1) with inputs the boundary values of the state $u(t, 0), u(t, 1)$ and the distributed perturbation term $f$ that appears in the PDE (2.1). An ISS Lyapunov functional for the PDE system (2.1) satisfies the following properties (see also [15] Chapter 1 for a slightly more demanding notion of an ISS Lyapunov functional for (2.1)):

**(i)** there exist functions $a_1, a_2 \in K_\infty$ such that the inequality $a_1(\|u\|_\infty) \leq V(u) \leq a_2(\|u\|_\infty)$ holds for all $u \in C^0([0,1])$,

**(ii)** for every $t_0 > 0$, $T > t_0$, the mapping $[t_0, T] \ni t \to V(u[t])$ is absolutely continuous for every solution $u[t]$ of (2.1),

**(iii)** there exist a set of functions $B \subseteq L^\infty(0,1)$, a semi-norm $G$ of $D := \Re \times \Re \times B$, a positive definite function $\rho \in C^0(\Re_+; \Re_+)$ and a function $\gamma \in K$ such that the following implication holds with $d[t] = (u(t,0), u(t,1), f[t]) \in D$ for every solution $u[t]$ of (2.1), for every $T > 0$, $t \in (0,T]$ for which $\frac{d}{dt} V(u[t])$ exists and every $\varepsilon > 0$: "if $(u[t], d[t]) \in C^0([0,1]) \times D$ satisfy $V(u[t]) \geq \gamma(G(d[t])) + \varepsilon$ then $\frac{d}{dt} V(u[t]) \leq -\rho(V(u[t]))$".

Property (i) holds for $V(u) = \|u\|_{\infty,\eta}$: notice that (2.5) implies that $\frac{\|u\|_\infty}{\max_{0 \leq x \leq 1}(\eta(x))} \leq \|u\|_{\infty,\eta} \leq \frac{\|u\|_\infty}{\min_{0 \leq x \leq 1}(\eta(x))}$

for all $u \in C^0([0,1])$. This explains why estimate (2.4) can give estimates of the sup norm of the state and this feature justifies the fact that the functional $V(u) = \|u\|_{\infty,\eta}$ is a coercive functional. The above properties enable us to derive estimate (2.4) by invoking the continuity of the mapping $[0,T] \ni t \to V(u[t])$. When $\zeta = 0$ then we get from (2.4) the estimate:

$$\|u[t]\|_{\infty,\eta} \leq \max\left(\|u[0]\|_{\infty,\eta}, \sup_{0 < s \leq t}\left(\max\left(\frac{|u(s,0)|}{\eta(0)}, \frac{|u(s,1)|}{\eta(1)}, \frac{\|f[s]\|_{\infty,\eta}}{\sigma}\right)\right)\right)$$

which is nothing else but a maximum principle for the solution of (2.1) (see [3,9,22]).

Theorem 2.1 holds for classical solutions of the PDE (2.1). Existence and uniqueness results for classical solutions of (nonlinear) parabolic PDEs are given in [3,17].

A disadvantage of Theorem 2.1 is the fact that Theorem 2.1 can provide ISS estimates only when boundary disturbances appear in Dirichlet boundary conditions. However, boundary disturbances can appear in other types of boundary conditions as well (e.g., Neumann or Robin boundary conditions). Our next main result provides an estimate of a weighted sup norm of the solution of the PDE (2.1) that depends on the initial condition, the distributed perturbation term that may be present in the PDE and the boundary values of the solution and its spatial derivative.

**Theorem 2.2:** Let $T > 0$, $a, b, c : (0,T) \times (0,1) \to \Re$, $f : (0,T] \times (0,1) \to \Re$ with $\sup_{0 < x < 1, 0 < t \leq T}(|f(t,x)|) < +\infty$, $u \in C^0([0,T] \times [0,1])$ with $u[t] \in C^2((0,1))$ for almost all $t \in (0,T)$, for which the derivative $\frac{\partial u}{\partial t}(t,x)$ exists and is a continuous function on $(0,T) \times [0,1]$, the derivatives $\frac{\partial u}{\partial x}(t,0)$, $\frac{\partial u}{\partial x}(t,1)$ exist for all $t \in (0,T]$ and equation (2.1) holds. Suppose that (2.2) holds and that there exist a constant $\sigma > 0$ and a positive function $\eta \in C^2([0,1];(0,+\infty))$ such that (2.3) holds. Then the following estimate holds for all $\zeta \in [0,\sigma)$, $g_0, g_1 : (0,T) \to (0,+\infty)$, $k_0, k_1 : (0,T) \to [1,+\infty)$ and $t \in (0,T]$:

$$\|u[t]\|_{\infty,\eta} \leq \max\left(\exp(-\zeta t)\|u[0]\|_{\infty,\eta}, \sup_{0 < s \leq t}\left(\max\left(r_0(s), r_1(s), \frac{\|f[s]\|_{\infty,\eta}}{\sigma - \zeta}\right)\exp(-\zeta(t-s))\right)\right) \quad (2.6)$$

where $\|u[t]\|_{\infty,\eta}$, $\|f[t]\|_{\infty,\eta}$ are defined by (2.5) and

$$r_0(t) := \min\left(\frac{|u(t,0)|}{\eta(0)}, \frac{g_0(t)}{\eta(0)}\left|\frac{\partial u}{\partial x}(t,0) - \left(\frac{\eta'(0)}{\eta(0)} + \frac{k_0(t)}{g_0(t)}\right)u(t,0)\right|\right), \text{ for } t \in (0,T] \quad (2.7)$$

$$r_1(t) := \min\left(\frac{|u(t,1)|}{\eta(1)}, \frac{g_1(t)}{\eta(1)}\left|\frac{\partial u}{\partial x}(t,1) - \left(\frac{\eta'(1)}{\eta(1)} - \frac{k_1(t)}{g_1(t)}\right)u(t,1)\right|\right), \text{ for } t \in (0,T] \quad (2.8)$$

It should be noted that estimate (2.6) is a more accurate ISS-style maximum principle estimate than estimate (2.4): notice that definitions (2.7), (2.8) imply that $r_0(t) \leq \frac{|u(t,0)|}{\eta(0)}$, $r_1(t) \leq \frac{|u(t,1)|}{\eta(1)}$ for all $t \in (0,T]$. However, there is a price to pay for obtaining the more accurate estimate (2.6): the solution has to be differentiable with respect to $x$ up to the boundary of the domain. This additional regularity requirement does not appear in Theorem 2.1. This additional regularity requirement is justified by the fact that Theorem 2.2 is a tool for handling boundary conditions which involve the spatial derivative at the boundary of the domain.

The functions $g_0, g_1 : (0,T] \to (0,+\infty)$, $k_0, k_1 : (0,T] \to [1,+\infty)$ are arbitrary and can be selected in an appropriate way in order to handle even nonlinear boundary conditions (see Example 3.2 in next section). However, when standard Neumann or Robin boundary conditions hold at the boundary then we obtain as a corollary of Theorem 2.2 the following result.

**Corollary 2.3:** *Let $T > 0$, $a,b,c : (0,T) \times (0,1) \to \Re$, $f : (0,T] \times (0,1) \to \Re$ with $\sup_{0<x<1, 0<t\leq T}(|f(t,x)|) < +\infty$, $u \in C^0([0,T] \times [0,1])$ with $u[t] \in C^2((0,1))$ for almost all $t \in (0,T)$, for which the derivative $\frac{\partial u}{\partial t}(t,x)$ exists and is a continuous function on $(0,T) \times [0,1]$, the derivatives $\frac{\partial u}{\partial x}(t,0)$, $\frac{\partial u}{\partial x}(t,1)$ exist for all $t \in (0,T]$ and equation (2.1) holds. Let $\mu_0, \mu_1 > 0$, $\lambda_0, \lambda_1 \in \Re$ be given constants. Suppose that (2.2) holds and that there exist a constant $\sigma > 0$ and a positive function $\eta \in C^2([0,1]; (0,+\infty))$ such that (2.3) holds. Moreover, suppose that one of the following conditions hold:*

$$\mu_0 \eta'(0) - \lambda_0 \eta(0) < 0 \quad (2.9)$$

$$\mu_1 \eta'(1) + \lambda_1 \eta(1) > 0 \quad (2.10)$$

$$\mu_0 \eta'(0) - \lambda_0 \eta(0) < 0 \text{ and } \mu_1 \eta'(1) + \lambda_1 \eta(1) > 0 \quad (2.11)$$

*Then estimate (2.6) holds for all $\zeta \in [0,\sigma)$ and $t \in (0,T]$, where*

$$r_0(t) := \min\left(\frac{|u(t,0)|}{\eta(0)}, \frac{1}{|\mu_0 \eta'(0) - \lambda_0 \eta(0)|}\left|\mu_0 \frac{\partial u}{\partial x}(t,0) - \lambda_0 u(t,0)\right|\right), \quad r_1(t) := \frac{|u(t,1)|}{\eta(1)}$$
$$\text{for } t \in (0,T], \text{ if (2.9) holds} \quad (2.12)$$

$$r_0(t) := \frac{|u(t,0)|}{\eta(0)}, \quad r_1(t) := \min\left(\frac{|u(t,1)|}{\eta(1)}, \frac{1}{\mu_1 \eta'(1) + \lambda_1 \eta(1)}\left|\mu_1 \frac{\partial u}{\partial x}(t,1) + \lambda_1 u(t,1)\right|\right),$$
$$\text{for } t \in (0,T], \text{ if (2.10) holds} \tag{2.13}$$

$$r_0(t) := \min\left(\frac{|u(t,0)|}{\eta(0)}, \frac{1}{|\mu_0 \eta'(0) - \lambda_0 \eta(0)|}\left|\mu_0 \frac{\partial u}{\partial x}(t,0) - \lambda_0 u(t,0)\right|\right),$$
$$r_1(t) := \min\left(\frac{|u(t,1)|}{\eta(1)}, \frac{1}{\mu_1 \eta'(1) + \lambda_1 \eta(1)}\left|\mu_1 \frac{\partial u}{\partial x}(t,1) + \lambda_1 u(t,1)\right|\right),$$
$$\text{for } t \in (0,T], \text{ if (2.11) holds} \tag{2.14}$$

The ISS-style maximum principle estimate (2.6), with $r_0, r_1$ defined by (2.12) or (2.13) or (2.14), provides in a direct way a fading memory ISS estimate when we have the boundary conditions

$$\mu_0 \frac{\partial u}{\partial x}(t,0) - \lambda_0 u(t,0) = d_0(t), \ u(t,1) = d_1(t), \text{ for all } t \in [0,T] \text{ (when (2.12) holds)}$$

or

$$u(t,0) = d_0(t), \ \mu_1 \frac{\partial u}{\partial x}(t,1) + \lambda_1 u(t,1) = d_1(t), \text{ for all } t \in [0,T] \text{ (when (2.13) holds)}$$

or

$$\mu_0 \frac{\partial u}{\partial x}(t,0) - \lambda_0 u(t,0) = d_0(t), \ \mu_1 \frac{\partial u}{\partial x}(t,1) + \lambda_1 u(t,1) = d_1(t), \text{ for all } t \in [0,T] \text{ (when (2.14) holds)}$$

where $d_0, d_1 \in C^0([0,T])$ are boundary disturbance inputs. However, it should be noted here that the ISS-style maximum principle estimate (2.6), with $r_0, r_1$ defined by (2.12) or (2.13) or (2.14), may also be useful in other cases which are rarely encountered in the literature. For example, if (2.11) holds and if the following highly nonlinear boundary conditions are valid

$$\left(\mu_0 \frac{\partial u}{\partial x}(t,0) - \lambda_0 u(t,0)\right) u(t,0) = 0, \ \left(\mu_1 \frac{\partial u}{\partial x}(t,1) + \lambda_1 u(t,1)\right) u(t,1) = 0, \text{ for all } t \in [0,T]$$

then the ISS-style maximum principle estimate (2.6), with $r_0, r_1$ defined by (2.14) gives the exponential decay estimate

$$\|u[t]\|_{\infty,\eta} \leq \exp(-\zeta t)\|u[0]\|_{\infty,\eta}$$

for zero source terms (i.e., when $f(t,x) \equiv 0$). This is a case where we can determine the asymptotic behavior of the solution even when we do not know exactly what happens at the boundary (i.e., for some $t \in [0,T]$ we may have $\mu_0 \frac{\partial u}{\partial x}(t,0) = \lambda_0 u(t,0)$ or $u(t,0) = 0$; similarly for $x=1$).

Conditions (2.9), (2.10), (2.11) are not new: they have also appeared in [11,15], where a completely different proof methodology was applied (the solution was approximated by a numerical scheme). However, here we do not follow the proof methodology in [11,15] and we obtain Corollary 2.3 under weaker requirements for the solution.

# 3. Illustrative Examples

In this section we present three examples that illustrate the use of the main results for the derivation of ISS estimates for nonlinear 1-D parabolic PDEs.

**Example 3.1 (A Nonlinear Reaction-Diffusion PDE):** Consider the nonlinear reaction-diffusion PDE

$$\frac{\partial u}{\partial t}(t,x) = \kappa(u(t,x))\frac{\partial^2 u}{\partial x^2}(t,x) + r(u(t,x))u(t,x) + f(t,x), \text{ for all } (t,x) \in (0,T]\times[0,1] \quad (3.1)$$

$$u(t,0) = d_0(t), \ u(t,1) = d_1(t), \text{ for all } t \in [0,T] \quad (3.2)$$

where $d_0, d_1 \in C^0([0,T])$ are boundary inputs, $f \in C^0([0,T]\times[0,1])$ is a distributed in-domain input and $\kappa, r \in C^0(\Re)$ are functions which satisfy the following conditions

$$\kappa(u) > 0, \text{ for all } u \in \Re \quad (3.3)$$

$$\sup\left\{\frac{\sigma + r(u)}{\kappa(u)} : u \in \Re\right\} < \pi^2 \quad (3.4)$$

where $\sigma > 0$ is a constant. We show next that there exists sufficiently small $\varphi > 0$ such that the following estimate holds for all $\zeta \in [0,\sigma)$ and $t \in (0,T]$

$$\|u[t]\|_{\infty,\eta} \leq \max\left(\exp(-\zeta t)\|u[0]\|_{\infty,\eta}, \sup_{0<s\leq t}\left(\max\left(\frac{|d_0(s)|}{\eta(0)}, \frac{|d_1(s)|}{\eta(1)}, \frac{\|f[s]\|_{\infty,\eta}}{\sigma - \zeta}\right)\exp(-\zeta(t-s))\right)\right) \quad (3.5)$$

with $\eta(x) = \sin(\theta x + \varphi)$, $x \in [0,1]$, for every $d_0, d_1 \in C^0([0,T])$, $f \in C^0([0,T]\times[0,1])$ and for every (classical) solution $u \in C^0([0,T]\times[0,1]) \cap C^1((0,T]\times[0,1])$ with $u[t] \in C^2([0,1])$ for $t \in (0,T]$ of (3.1), (3.2).

Indeed, inequality (3.4) shows that there exist a constant $\theta \in (0,\pi)$ (sufficiently close to $\pi$) such that

$$\theta^2 \geq \frac{\sigma + r(u)}{\kappa(u)}, \text{ for all } u \in \Re \quad (3.6)$$

Let $\varphi > 0$ be a constant sufficiently small so that $\varphi + \theta < \pi$. Notice that for every solution of (3.1), (3.2), equation (2.1) holds with

$$a(t,x) = \kappa(u(t,x)), \ b(t,x) = 0, \ c(t,x) = r(u(t,x)). \quad (3.7)$$

Inequalities (3.3), (3.6) and definitions (3.7) guarantee that (2.2) and (2.3) hold with $\eta(x) = \sin(\theta x + \varphi)$. Theorem 2.1 guarantees that for every $T > 0$, for every $d_0, d_1 \in C^0([0,T])$, $f \in C^0([0,T]\times[0,1])$ and for every (classical) solution $u \in C^0([0,T]\times[0,1]) \cap C^1((0,T]\times[0,1])$ with $u[t] \in C^2([0,1])$ for $t \in (0,T]$ of the PDE problem (3.1), (3.2), estimate (3.5) holds.

It should be noted here that the stability conditions (3.3), (3.4) give a sharp characterization of the stability for the PDE problem (3.1), (3.2). Indeed, when $\kappa(u) \equiv p > 0$, $r(u) \equiv p\pi^2$ (and consequently $\sup\left\{\frac{r(u)}{\kappa(u)} : u \in \Re\right\} = \pi^2$), $d_0 = d_1 \equiv 0$, $f \equiv 0$ then there are solutions of the (linear) PDE problem (3.1), (3.2) that do not tend to zero (namely solutions of the form $u(t,x) = A\sin(\pi x)$ where $A \in \Re$). Therefore, estimate (3.5) cannot hold in this case. ◁

The following example illustrates how Theorem 2.2 can be used for the derivation of ISS estimates for parabolic PDEs with nonlinear and non-local boundary conditions. Moreover, the example shows that we can study parabolic PDEs for which the diffusion coefficient is a non-local functional of the state. This situation is important because it arises in many phenomena related to turbulent flows (see [4,5]) as well as in the work of O. A. Ladyzhenskaya (see [18]) on the modification of the Navier-Stokes equations and in feedback control of fluid flows (see [1]).

**Example 3.2 (Heat Equation with Nonlinear and Non-Local Boundary Conditions):** Consider the following nonlinear and non-local 1-D parabolic PDE:

$$\frac{\partial u}{\partial t}(t,x) = \kappa(u[t])\frac{\partial^2 u}{\partial x^2}(t,x) + f(t,x), \text{ for all } (t,x) \in (0,T] \times [0,1] \tag{3.8}$$

$$\frac{\partial u}{\partial x}(t,0) = \left(\lambda_0 + \beta_0(u[t])\right)u(t,0) + d_0(t), \quad \frac{\partial u}{\partial x}(t,1) = -\left(\lambda_1 + \beta_1(u[t])\right)u(t,1) + d_1(t) \tag{3.9}$$

where $d_0, d_1 \in C^0([0,T])$ are the boundary disturbances, $f \in C^0([0,T] \times [0,1])$ is the distributed disturbance, $\kappa : C^1([0,1]) \to \Re$, $\beta_0, \beta_1 : C^0([0,1]) \to \Re_+$ are given non-negative functionals and $\lambda_0, \lambda_1 > 0$ are constants.

We assume that there exists a constant $\kappa^* > 0$ for which $\kappa(u) \geq \kappa^*$ for all $u \in C^0([0,1])$. Let $\theta \in \left(0, \frac{\pi}{2}\right)$ be a constant sufficiently small so that $\theta\tan(\theta) < \lambda_1$. We show next that every classical solution of (3.8), (3.9) $u \in C^0([0,T] \times [0,1]) \cap C^1((0,T] \times [0,1])$ with $u[t] \in C^2([0,1])$ for $t \in (0,T]$ satisfies the following estimate holds for all $\zeta \in [0, \kappa^*\theta^2)$ and $t \in (0,T]$

$$\|u[t]\|_{\infty,\eta}$$

$$\leq \max\left(\exp(-\zeta t)\|u[0]\|_{\infty,\eta}, \sup_{0 < s \leq t}\left(\max\left(\frac{|d_0(s)|}{\lambda_0}, \frac{|d_1(s)|}{\lambda_1\cos(\theta) - \theta\sin(\theta)}, \frac{\|f[s]\|_{\infty,\eta}}{\kappa^*\theta^2 - \zeta}\right)\exp(-\zeta(t-s))\right)\right)$$

(3.10)

with $\eta(x) = \cos(\theta x)$ for $x \in [0,1]$. Indeed, notice that for every solution of (3.8), (3.9), equation (2.1) holds with

$$a(t,x) = \kappa(u[t]), \quad b(t,x) \equiv 0, \quad c(t,x) \equiv 0. \tag{3.11}$$

Moreover, due to the fact that $\kappa(u) \geq \kappa^*$ for all $u \in C^1([0,1])$, inequalities (2.2) and (2.3) hold with $\sigma = \kappa^*\theta^2$ and $\eta(x) = \cos(\theta x)$. Estimate (3.10) is a direct consequence of Theorem 2.2 with $g_0(t) = \frac{1}{\beta_0(u[t]) + \lambda_0}$, $k_0(t) = k_1(t) \equiv 1$ and $g_1(t) = \frac{1}{\beta_1(u[t]) + \lambda_1 - \theta\tan(\theta)}$. ◁

The final example shows that the main results of the present work can be used for the study of realistic problems where the physical properties of a material are nonlinear functions of thermodynamic variables (such as temperature or pressure).

**Example 3.3 (The "real" heat equation):** Thermal conductivity $k$ and volumetric heat capacity $c$ are quantities that depend on the temperature $u$ of a solid material. Consequently, by exploiting the conservation of energy and Fourier's law, we obtain the following equation for a 1-D spatial domain:

$$\frac{\partial u}{\partial t}(t,x) = \frac{k(u(t,x))}{c(u(t,x))} \frac{\partial^2 u}{\partial x^2}(t,x) + \frac{k'(u(t,x))}{c(u(t,x))} \left( \frac{\partial u}{\partial x}(t,x) \right)^2 , \text{ for all } (t,x) \in (0,T] \times [0,1] \quad (3.12)$$

Equation (3.12) is important for certain solids when the temperature variation is large. Therefore, we are led to the study of the PDE problem

$$\frac{\partial u}{\partial t}(t,x) = \kappa(u(t,x)) \frac{\partial^2 u}{\partial x^2}(t,x) + g(u(t,x)) \left( \frac{\partial u}{\partial x}(t,x) \right)^2 , \text{ for all } (t,x) \in (0,T] \times [0,1] \quad (3.13)$$

$$u(t,0) = d_0(t) , \ u(t,1) = d_1(t) , \text{ for all } t \in [0,T] \quad (3.14)$$

where $d_0, d_1 \in C^0([0,T])$ are boundary disturbances and $\kappa, g \in C^0(\Re)$ are given functions with $\kappa(u) > 0$ for all $u \in \Re$. We assume that the increasing function $\gamma : \Re \to \Re$ defined by

$$\gamma(u) := \int_0^u \exp\left( \int_0^s \frac{g(l)}{\kappa(l)} dl \right) ds , \text{ for } u \in \Re \quad (3.15)$$

satisfies the condition $\gamma(\Re) = \Re$ or equivalently,

$$\lim_{u \to \pm\infty} (\gamma(u)) = \pm\infty \quad (3.16)$$

Moreover, we assume that there exists $\kappa^* > 0$ such that

$$\kappa(u) \geq \kappa^* \text{ for all } u \in \Re . \quad (3.17)$$

It should be noticed here that the "real" heat equation problem (3.12) is accompanied by constraints of the form $u \in (u_{\min}, u_{\max})$, where $u_{\min}$ corresponds to the absolute zero temperature (3$^{\text{rd}}$ law of thermodynamics) and $u_{\max}$ corresponds to the melting point. In this case, one can smoothly extend the thermal conductivity $k$ and volumetric heat capacity $c$ out of the interval $(u_{\min}, u_{\max})$ so that conditions (3.16), (3.17) hold with $\kappa(u) = k(u)/c(u)$ and $g(u) = k'(u)/c(u)$.

The invertible transformation $w(t,x) = \gamma(u(t,x))$ transforms the PDE problem (3.13), (3.14) to the following PDE problem

$$\frac{\partial w}{\partial t}(t,x) = \kappa\left(\gamma^{-1}(w(t,x))\right) \frac{\partial^2 w}{\partial x^2}(t,x) , \text{ for all } (t,x) \in (0,T] \times [0,1] \quad (3.18)$$

$$w(t,0) = \tilde{d}_0(t) = \gamma(d_0(t)), \ w(t,1) = \tilde{d}_1(t) = \gamma(d_1(t)), \text{ for all } t \in [0,T] \quad (3.19)$$

Following the analysis presented in Example 3.1, for every $\varphi \in (0, \pi/2)$ every (classical) solution $w \in C^0([0,T] \times [0,1]) \cap C^1((0,T] \times [0,1])$ with $w[t] \in C^2([0,1])$ for $t \in (0,T]$ of (3.18), (3.19) satisfies the following estimate holds for all $\zeta \in \left(0, \kappa^*(\pi - 2\varphi)^2\right)$ and $t \in (0,T]$

$$\|w[t]\|_{\infty,\eta} \leq \max\left(\exp(-\zeta t)\|w[0]\|_{\infty,\eta}, \frac{1}{\sin(\varphi)} \sup_{0<s\leq t}\left(\max\left(|\tilde{d}_0(s)|,|\tilde{d}_1(s)|\right)\exp(-\zeta(t-s))\right)\right) \quad (3.20)$$

with $\eta(x) = \sin((\pi - 2\varphi)x + \varphi)$ for $x \in [0,1]$. Exploiting definition (2.5), we get from (3.20) for all $\zeta \in (0, \kappa^*(\pi - 2\varphi)^2)$ and $t \in (0,T]$:

$$\|w[t]\|_{\infty} \leq \frac{1}{\sin(\varphi)} \max\left(\exp(-\zeta t)\|w[0]\|_{\infty}, \sup_{0<s\leq t}\left(\max\left(|\tilde{d}_0(s)|,|\tilde{d}_1(s)|\right)\exp(-\zeta(t-s))\right)\right) \quad (3.21)$$

Next define:

$$\gamma_1(s) := \min(\gamma(s), -\gamma(-s)), \text{ for } s \geq 0 \quad (3.22)$$

$$\gamma_2(s) := \max(\gamma(s), -\gamma(-s)), \text{ for } s \geq 0 \quad (3.23)$$

Condition (3.16) guarantees that both functions $\gamma_1, \gamma_2$ are of class $K_\infty$ and satisfy the following property:

$$\gamma_1(|u|) \leq |\gamma(u)| \leq \gamma_2(|u|), \text{ for all } u \in \Re \quad (3.24)$$

Using the fact that $w(t,x) = \gamma(u(t,x))$, (3.19) and (3.24), we obtain from (3.21) for all $\zeta \in (0, \kappa^*(\pi - 2\varphi)^2)$ and $t \in (0,T]$:

$$\|u[t]\|_{\infty} \leq \max\left(\omega_\zeta(\|u[0]\|_{\infty}, t), \sup_{0<s\leq t}\left(\max\left(\omega_\zeta(|d_0(s)|, t-s), \omega_\zeta(|d_1(s)|, t-s)\right)\right)\right) \quad (3.25)$$

where $\omega_\zeta(s,t) := \gamma_1^{-1}\left(\frac{\exp(-\zeta t)}{\sin(\varphi)}\gamma_2(s)\right)$. The fading memory ISS estimate (3.25) shows that the gain functions for each one of the boundary disturbances is $\omega_\zeta(s,0) := \gamma_1^{-1}\left(\frac{\gamma_2(s)}{\sin(\varphi)}\right)$. Crystalline semiconductors, crystalline quartz and crystalline ceramics, used in the manufacture of quartz clocks and various kinds of electronic devices, including diodes, transistors, and integrated circuits, are known to have a thermal conductivity that scales as $1/u$ (see [23]). With further calculations, which would need to be mindful of the absolute zero and melting point limits on the temperature, one could pursue deriving a growth rate on the ISS gain $\omega_\zeta(s,0) := \gamma_1^{-1}\left(\frac{\gamma_2(s)}{\sin(\varphi)}\right)$. ◁

## 4. Proofs of Main Results

For the proofs of the main results of the paper we need some auxiliary results. The first auxiliary result provides a sufficient condition that guarantees that the sup norm is an absolutely continuous function of time.

**Lemma 4.1:** *Let $b > a$ be given real numbers and let $u \in C^0([a,b] \times [0,1])$ be a given function for which the derivative $\frac{\partial u}{\partial t}(t,x)$ exists and is a continuous function on $[a,b] \times [0,1]$. Then the mapping $[a,b] \ni t \to \|u[t]\|_{\infty}$ is a Lipschitz mapping. Moreover, for every $t \in [a,b]$ for which the limit $\lim_{h \to 0^+}\left(\frac{\|u[t+h]\|_{\infty} - \|u[t]\|_{\infty}}{h}\right)$ exists, the following equality holds:*

$$\lim_{h \to 0^+} \left( \frac{\|u[t+h]\|_\infty - \|u[t]\|_\infty}{h} \right) = \lim_{h \to 0^+} \left( \frac{\left\| u[t] + h \frac{\partial u}{\partial t}[t] \right\|_\infty - \|u[t]\|_\infty}{h} \right) \quad (4.1)$$

**Proof:** For every $t_1, t_2 \in [a,b]$ with $t_1 \leq t_2$ we get:

$$\|u[t_2]\|_\infty = \max_{0 \leq x \leq 1} \left( |u(t_2, x)| \right) = \max_{0 \leq x \leq 1} \left( \left| u(t_1, x) + \int_{t_1}^{t_2} \frac{\partial u}{\partial t}(s, x) ds \right| \right)$$

$$\leq \max_{0 \leq x \leq 1} \left( |u(t_1, x)| + \int_{t_1}^{t_2} \left| \frac{\partial u}{\partial t}(s, x) \right| ds \right)$$

$$\leq \max_{0 \leq x \leq 1} \left( |u(t_1, x)| + (t_2 - t_1) \max_{a \leq s \leq b, 0 \leq x \leq 1} \left( \left| \frac{\partial u}{\partial t}(s, x) \right| \right) \right)$$

$$= \|u[t_1]\|_\infty + (t_2 - t_1) \max_{a \leq s \leq b, 0 \leq x \leq 1} \left( \left| \frac{\partial u}{\partial t}(s, x) \right| \right)$$

Similarly as above, we get:

$$\|u[t_1]\|_\infty = \max_{0 \leq x \leq 1} \left( |u(t_1, x)| \right) = \max_{0 \leq x \leq 1} \left( \left| u(t_2, x) - \int_{t_1}^{t_2} \frac{\partial u}{\partial t}(s, x) ds \right| \right)$$

$$\leq \max_{0 \leq x \leq 1} \left( |u(t_2, x)| + \int_{t_1}^{t_2} \left| \frac{\partial u}{\partial t}(s, x) \right| ds \right)$$

$$\leq \max_{0 \leq x \leq 1} \left( |u(t_2, x)| + (t_2 - t_1) \max_{a \leq s \leq b, 0 \leq x \leq 1} \left( \left| \frac{\partial u}{\partial t}(s, x) \right| \right) \right)$$

$$= \|u[t_2]\|_\infty + (t_2 - t_1) \max_{a \leq s \leq b, 0 \leq x \leq 1} \left( \left| \frac{\partial u}{\partial t}(s, x) \right| \right)$$

Consequently, the two above inequalities guarantee that the following inequality holds for all $t_1, t_2 \in [a,b]$:

$$\left| \|u[t_2]\|_\infty - \|u[t_1]\|_\infty \right| \leq |t_2 - t_1| \max_{a \leq s \leq b, 0 \leq x \leq 1} \left( \left| \frac{\partial u}{\partial t}(s, x) \right| \right) \quad (4.2)$$

Inequality (4.2) shows that the mapping $[a,b] \ni t \to \|u[t]\|_\infty$ is a Lipschitz mapping.

Let $t \in [a,b)$, $\varepsilon > 0$ be given (arbitrary). By virtue of uniform continuity of $\frac{\partial u}{\partial t}(t, x)$ on $[a,b] \times [0,1]$ (which follows from the compactness of $[a,b] \times [0,1]$) there exists $h > 0$ sufficiently small such that $\left| \frac{\partial u}{\partial t}(s, x) - \frac{\partial u}{\partial t}(t, x) \right| \leq \varepsilon$ for all $x \in [0,1]$, $s, t \in [a,b]$ with $|s - t| \leq h$. Let $t \in [a,b)$ be given (arbitrary). We get:

$$\|u[t+h]\|_\infty = \max_{0\le x\le 1}\left(\left|u(t,x)+h\frac{\partial u}{\partial t}(t,x)+\int_t^{t+h}\left(\frac{\partial u}{\partial t}(s,x)-\frac{\partial u}{\partial t}(t,x)\right)ds\right|\right)$$

$$\le \left\|u[t]+h\frac{\partial u}{\partial t}[t]\right\|_\infty + \varepsilon h$$

Similarly as above, we get:

$$\|u[t+h]\|_\infty = \max_{0\le x\le 1}\left(\left|u(t,x)+h\frac{\partial u}{\partial t}(t,x)+\int_t^{t+h}\left(\frac{\partial u}{\partial t}(s,x)-\frac{\partial u}{\partial t}(t,x)\right)ds\right|\right)$$

$$\ge \max_{0\le x\le 1}\left(\left|u(t,x)+h\frac{\partial u}{\partial t}(t,x)\right|-\int_t^{t+h}\left|\frac{\partial u}{\partial t}(s,x)-\frac{\partial u}{\partial t}(t,x)\right|ds\right)$$

$$\ge \left\|u[t]+h\frac{\partial u}{\partial t}[t]\right\|_\infty - \varepsilon h$$

Since $\varepsilon > 0$ is arbitrary, the above inequalities allow us to conclude that

$$\lim_{h\to 0^+}\left(\frac{\left\|u[t]+h\frac{\partial u}{\partial t}[t]\right\|_\infty - \|u[t+h]\|_\infty}{h}\right) = 0, \text{ for all } t\in [a,b) \quad (4.3)$$

Equality (4.1) is a direct consequence of equation (4.3). ◁

The second auxiliary result is a technical lemma that provides a formula for a specific limit. This specific limit was encountered in the statement of Lemma 4.1 and consequently, the formula that the following lemma provides, plays a crucial role.

**Lemma 4.2:** *Let $u, w \in C^0([0,1])$ be given functions and suppose that the limit $\lim_{h\to 0^+}\left(\frac{\|u+hw\|_\infty - \|u\|_\infty}{h}\right)$ exists. Then $\lim_{h\to 0^+}\left(\frac{\|u+hw\|_\infty - \|u\|_\infty}{h}\right) \le \|w\|_\infty$. Moreover, if $\|u\|_\infty > 0$ then*

$$\lim_{h\to 0^+}\left(\frac{\|u+hw\|_\infty - \|u\|_\infty}{h}\right) \le \max_{x\in I}\left(\text{sgn}(u(x))w(x)\right) \quad (4.4)$$

*where $I = \left\{y\in [0,1]: |u(y)| = \|u\|_\infty\right\}$.*

**Proof:** The fact that $\lim_{h\to 0^+}\left(\frac{\|u+hw\|_\infty - \|u\|_\infty}{h}\right) \le \|w\|_\infty$ follows from the inequality $\|u+hw\|_\infty \le \|u\|_\infty + h\|w\|_\infty$ that holds for all $h > 0$.

Next we assume that $\|u\|_\infty > 0$ and we show that (4.4) holds. Let $\varepsilon \in \left(0, \|u\|_\infty\right)$ be given and define the sets:

$$I_\varepsilon^+ = \left\{y\in [0,1]: u(y) \ge \|u\|_\infty - \varepsilon\right\},\ I_\varepsilon^- = \left\{y\in [0,1]: u(y) \le -\|u\|_\infty + \varepsilon\right\},\ I_\varepsilon = I_\varepsilon^+ \cup I_\varepsilon^- \quad (4.5)$$

By virtue of continuity of $u$ each of the sets $I_\varepsilon, I_\varepsilon^+, I_\varepsilon^-$ is compact and $I_\varepsilon$ is non-empty (but notice that one of the sets $I_\varepsilon^+, I_\varepsilon^-$ may be empty). Definitions (4.5) imply that for all $y \in [0,1] \setminus I_\varepsilon$ we have:

$$|u(y) + hw(y)| \leq |u(y)| + h|w(y)| \leq \|u\|_\infty - \varepsilon + h\|w\|_\infty \tag{4.6}$$

Since $\|u + hw\|_\infty \geq \|u\|_\infty - h\|w\|_\infty$, we conclude that for $y \in [0,1] \setminus I_\varepsilon$ and $h > 0$ sufficiently small (so that $2h\|w\|_\infty < \varepsilon$), it holds that $|u(y) + hw(y)| < \|u + hw\|_\infty$. Thus we conclude that $\|u + hw\|_\infty = \max_{x \in I_\varepsilon}(|u(x) + hw(x)|)$ for all $h > 0$ sufficiently small, which also implies that

$$\begin{aligned}\|u + hw\|_\infty &= \max_{x \in I_\varepsilon}(|u(x) + hw(x)|) = \max_{x \in I_\varepsilon}(|u(x) + hw(x)| - |u(x)| + |u(x)|) \\ &\leq \max_{x \in I_\varepsilon}(|u(x) + hw(x)| - |u(x)|) + \|u\|_\infty \\ &= h\max_{x \in I_\varepsilon}(\mathrm{sgn}(u(x))w(x)) + \|u\|_\infty\end{aligned} \tag{4.7}$$

for all $h > 0$ sufficiently small. The last equality in (4.7) is obtained by distinguishing the cases $x \in I_\varepsilon^+$, $x \in I_\varepsilon^-$, using definitions (4.5) and assuming that $h\|w\|_\infty < \|u\|_\infty - \varepsilon$. Inequality (4.7) and the fact that $\varepsilon \in (0, \|u\|_\infty)$ is arbitrary imply the following inequality:

$$\lim_{h \to 0^+}\left(\frac{\|u + hw\|_\infty - \|u\|_\infty}{h}\right) \leq \max_{x \in I_\varepsilon}(\mathrm{sgn}(u(x))w(x)), \text{ for all } \varepsilon \in (0, \|u\|_\infty) \tag{4.8}$$

We show next by contradiction that (4.4) holds. Assume that $\lim_{h \to 0^+}\left(\frac{\|u + hw\|_\infty - \|u\|_\infty}{h}\right) > \max_{x \in I}(\mathrm{sgn}(u(x))w(x))$. Define

$$p = \lim_{h \to 0^+}\left(\frac{\|u + hw\|_\infty - \|u\|_\infty}{h}\right) - \max_{x \in I}(\mathrm{sgn}(u(x))w(x)) > 0 \tag{4.9}$$

Consider a sequence $\{y_n \in [0,1]\}_{n=0}^\infty$ that satisfies $y_n \in I_\varepsilon$ and $\mathrm{sgn}(u(y_n))w(y_n) = \max_{x \in I_\varepsilon}(\mathrm{sgn}(u(x))w(x))$ with $\varepsilon = 2^{-n-1}\|u\|_\infty$ for $n = 0,1,2,\ldots$. Clearly, this sequence is bounded and consequently, it contains a convergent subsequence. Without loss of generality, we denote this convergent subsequence by $\{x_n \in [0,1]\}_{n=0}^\infty \subseteq \{y_n \in [0,1]\}_{n=0}^\infty$. Define $\lim_{n \to +\infty}(x_n) = x^*$ and notice that since for each $n = 0,1,2,\ldots$ there exists $m \geq n$ with $x_n = y_m \in I_\varepsilon$ and $\varepsilon = 2^{-m-1}\|u\|_\infty$, it follows from definitions (4.5) that $\|u\|_\infty(1 - 2^{-n-1}) \leq \|u\|_\infty(1 - 2^{-m-1}) \leq |u(x_n)|$. Continuity of $u$ and the fact that $I = \{y \in [0,1] : |u(y)| = \|u\|_\infty\}$ imply that $x^* \in I$.

Since for each $n = 0,1,2,\ldots$ there exists $m \geq n$ with $x_n = y_m \in I_\varepsilon$ and $\mathrm{sgn}(u(y_m))w(y_m) = \max_{x \in I_\varepsilon}(\mathrm{sgn}(u(x))w(x))$ with $\varepsilon = 2^{-m-1}\|u\|_\infty$, it follows from (4.8) and (4.9) that

$$\mathrm{sgn}(u(x_n))w(x_n) \geq p + \max_{x \in I}(\mathrm{sgn}(u(x))w(x)) \text{ for all } n = 0,1,2,\ldots. \tag{4.10}$$

Since $x^* \in I$ we must have either $u(x^*) = \|u\|_\infty > 0$ or $u(x^*) = -\|u\|_\infty < 0$. Consequently, since $\lim_{n \to +\infty}(x_n) = x^*$ and since $u$ is continuous, we must have either $u(x_n) > 0$ or $u(x_n) < 0$ for sufficiently large $n$. When $u(x_n) > 0$ for sufficiently large $n$, we get by virtue of continuity of $w$ and (4.10) that

$$\text{sgn}(u(x_n))w(x_n) \to w(x^*) = \text{sgn}(u(x^*))w(x^*) \geq p + \max_{x \in I}(\text{sgn}(u(x))w(x))$$

The above inequality together with the fact that $\text{sgn}(u(x^*))w(x^*) \leq \max_{x \in I}(\text{sgn}(u(x))w(x))$ (a consequence of the fact that $x^* \in I$) implies that $p \leq 0$; a contradiction with (9). Similarly, we obtain a contradiction when $u(x_n) < 0$ for sufficiently large $n$.

The proof is complete.  ◁

We are now ready to derive estimates for the solution of the PDE (2.1) under assumption (2.2). The third auxiliary result provides an estimate under an additional assumption for the coefficient of the reaction term $c(t,x)$.

**Lemma 4.3:** *Let $T > 0$, $a,b,c : (0,T) \times (0,1) \to \Re$, $f : (0,T] \times (0,1) \to \Re$ with $\sup_{0<x<1, 0<t\leq T}(|f(t,x)|) < +\infty$, $u \in C^0([0,T] \times [0,1])$ with $u[t] \in C^2((0,1))$ for almost all $t \in (0,T)$, for which the derivative $\frac{\partial u}{\partial t}(t,x)$ exists and is a continuous function on $(0,T) \times [0,1]$ and for which equation (2.1) holds. Suppose that inequality (2.2) holds and that*

$$\sigma := -\sup_{0<x<1, t\in(0,T)}(c(t,x)) > 0 \tag{4.11}$$

*Then the following estimate holds for all $\zeta \in [0,\sigma)$ and $t \in (0,T]$:*

$$\|u[t]\|_\infty \leq \max\left(\exp(-\zeta t)\|u[0]\|_\infty, \sup_{0<s\leq t}\left(\max\left(|u(s,0)|, |u(s,1)|, \frac{\|f[s]\|_\infty}{\sigma - \zeta}\right)\exp(-\zeta(t-s))\right)\right) \tag{4.12}$$

**Proof:** Let (arbitrary) $\zeta \in (0,\sigma)$, $t_1, t_2 \in (0,T)$ with $t_1 < t_2$ be given. Lemma 4.1 implies that the mapping $[t_1, t_2] \ni t \to \|u[t]\|_\infty$ is an absolutely continuous mapping. Moreover, there exists a Lebesgue measure zero set $N \subset [t_1, t_2]$ such that the following equation holds:

$$\frac{d}{dt}(\|u[t]\|_\infty) = \lim_{h \to 0^+}\left(\frac{\left\|u[t] + h\frac{\partial u}{\partial t}[t]\right\|_\infty - \|u[t]\|_\infty}{h}\right), \text{ for all } t \in [t_1, t_2] \setminus N \tag{4.13}$$

Moreover, we assume that the Lebesgue measure zero set $N \subset [t_1, t_2]$ has been selected in such a way that (2.1) holds for all $t \in [t_1, t_2] \setminus N$ and for all $x \in (0,1)$ and $u[t] \in C^2((0,1))$ for all $t \in [t_1, t_2] \setminus N$.

Lemma 4.2 implies that

$$\frac{d}{dt}\left(\|u[t]\|_\infty\right) \leq \max_{x \in I(t)}\left(\operatorname{sgn}(u(t,x))\frac{\partial u}{\partial t}(t,x)\right), \text{ for all } t \in [t_1,t_2]\setminus N \text{ with } \|u[t]\|_\infty > 0 \quad (4.14)$$

where $I(t) = \{ y \in [0,1] : |u(t,y)| = \|u[t]\|_\infty \}$.

Let (arbitrary) $\varepsilon > 0$ be given. Inequality (4.14) guarantees that the following implication holds for all $t \in [t_1, t_2]\setminus N$:

$$\|u[t]\|_\infty \geq \varepsilon + \max\left(|u(t,0)|, |u(t,1)|, \frac{\|f[t]\|_\infty}{\sigma - \zeta}\right) \Rightarrow \frac{d}{dt}\left(\|u[t]\|_\infty\right) \leq \max_{x \in I(t)}\left(\operatorname{sgn}(u(t,x))\frac{\partial u}{\partial t}(t,x)\right) \quad (4.15)$$

Pick any $t \in [t_1,t_2]\setminus N$. If $x \in I(t)$ is an interior point of $[0,1]$ (i.e., $x \in (0,1)$) and satisfies $u(t,x) = \|u[t]\|_\infty$ then $u$ has a maximum at $x$ and consequently (since $u[t] \in C^2((0,1))$) $\frac{\partial^2 u}{\partial x^2}(t,x) \leq 0$, $\frac{\partial u}{\partial x}(t,x) = 0$. It follows from (2.1) and (2.2) that $\operatorname{sgn}(u(t,x))\frac{\partial u}{\partial t}(t,x) \leq c(t,x)\|u[t]\|_\infty + f(t,x)$ for every $x \in I(t) \cap (0,1)$ with $u(t,x) = \|u[t]\|_\infty$. If $x \in I(t)$ is an interior point of $[0,1]$ (i.e., $x \in (0,1)$) and satisfies $u(t,x) = -\|u[t]\|_\infty$ then $u$ has a minimum at $x$ and consequently (since $u[t] \in C^2((0,1))$) $\frac{\partial^2 u}{\partial x^2}(t,x) \geq 0$, $\frac{\partial u}{\partial x}(t,x) = 0$. It follows from (2.1) and (2.2) that $\operatorname{sgn}(u(t,x))\frac{\partial u}{\partial t}(t,x) \leq c(t,x)\|u[t]\|_\infty - f(t,x)$ for every $x \in I(t) \cap (0,1)$ with $u(t,x) = -\|u[t]\|_\infty$. Combining both cases, we obtain:

If $t \in [t_1,t_2]\setminus N$ and $x \in I(t) \cap (0,1)$ then $\operatorname{sgn}(u(t,x))\frac{\partial u}{\partial t}(t,x) \leq c(t,x)\|u[t]\|_\infty + |f(t,x)|$ \quad (4.16)

Notice that the inequality $\|u[t]\|_\infty \geq \varepsilon + \max(|u(t,0)|, |u(t,1)|)$ implies that $I(t) \subseteq (0,1)$. Consequently, we obtain from (4.11), (4.15) and (4.16) for all $t \in [t_1,t_2]\setminus N$:

$$\|u[t]\|_\infty \geq \varepsilon + \max\left(|u(t,0)|, |u(t,1)|, \frac{\|f[t]\|_\infty}{\sigma - \zeta}\right) \Rightarrow \frac{d}{dt}\left(\|u[t]\|_\infty\right) \leq -\zeta\|u[t]\|_\infty \quad (4.17)$$

Implication (4.17) and Lemma 2.14 on page 82 in [10] imply that the following estimate holds for all $t \in [t_1, t_2]$:

$$\|u[t]\|_\infty \leq \max\left(\exp(-\zeta(t-t_1))\|u[t_1]\|_\infty, \sup_{t_1 \leq s \leq t}\left(\max\left(|u(s,0)|, |u(s,1)|, \frac{\|f[s]\|_\infty}{\sigma - \zeta}\right)\exp(-\zeta(t-s))\right)\right) + \varepsilon$$
(4.18)

Since $\varepsilon > 0$ is arbitrary, we conclude from (4.18) that the following estimate holds for all $t \in [t_1, t_2]$:

$$\|u[t]\|_\infty \leq \max\left(\exp(-\zeta(t-t_1))\|u[t_1]\|_\infty, \sup_{0 < s \leq t}\left(\max\left(|u(s,0)|, |u(s,1)|, \frac{\|f[s]\|_\infty}{\sigma - \zeta}\right)\exp(-\zeta(t-s))\right)\right)$$
(4.19)

Finally, since $t_1, t_2 \in (0,T)$ are arbitrary and since $u \in C^0([0,T]\times[0,1])$ (which implies that $\lim_{t_1 \to 0^+}(\|u[t_1]\|_\infty) = \|u[0]\|_\infty$ and $\lim_{t_2 \to T^-}(\|u[t_2]\|_\infty) = \|u[T]\|_\infty$), we conclude from (4.19) that (4.12) holds

for all $t \in (0,T]$ and $\zeta \in (0,\sigma)$. Continuity arguments guarantee the fact that estimate (4.12) holds for $\zeta = 0$ as well. The proof is complete. ◁

We are now ready to give the proof of Theorem 2.1.

**Proof of Theorem 2.1:** Theorem 2.1 is a direct consequence of Lemma 4.3 applied to the function $w \in C^0([0,T] \times [0,1])$ defined by

$$w(t,x) = u(t,x)/\eta(x), \text{ for all } t \in [0,T],\ x \in [0,1] \quad (4.20)$$

The function $w \in C^0([0,T] \times [0,1])$ satisfies the following properties:

- $w[t] \in C^2((0,1))$ for almost all $t \in (0,T)$,
- the derivative $\dfrac{\partial w}{\partial t}(t,x)$ exists, is a continuous function on $(0,T) \times [0,1]$,
- the following equation holds:

$$\frac{\partial w}{\partial t}(t,x) = a(t,x)\frac{\partial^2 w}{\partial x^2}(t,x) + \left(2a(t,x)\frac{\eta'(x)}{\eta(x)} + b(t,x)\right)\frac{\partial w}{\partial x}(t,x)$$
$$+ \frac{1}{\eta(x)}\big(a(t,x)\eta''(x) + b(t,x)\eta'(x) + c(t,x)\eta(x)\big)w(t,x) + \frac{f(t,x)}{\eta(x)}$$
$$\text{for almost all } t \in (0,T) \text{ and for all } x \in (0,1) \quad (4.21)$$

Notice that definition (4.20) implies that $\|u[t]\|_{\infty,\eta} = \|w[t]\|_\infty$ for all $t \in [0,T]$.

The proof is complete. ◁

For the proof of Theorem 2.2 we need an additional auxiliary result, which provides an estimate for the solution of (2.1) under (2.2) and (4.11). The estimate is different from the one provided by Lemma 4.3 and the solution of (2.1) is assumed to be differentiable with respect to $x$ at the boundary of the domain.

**Lemma 4.4:** *Let $T > 0$, $a,b,c : (0,T) \times (0,1) \to \Re$, $f : (0,T] \times (0,1) \to \Re$ with $\sup_{0<x<1, 0<t\leq T}(|f(t,x)|) < +\infty$, $u \in C^0([0,T] \times [0,1])$ with $u[t] \in C^2((0,1))$ for almost all $t \in (0,T)$, for which the derivative $\dfrac{\partial u}{\partial t}(t,x)$ exists and is a continuous function on $(0,T) \times [0,1]$, the derivatives $\dfrac{\partial u}{\partial x}(t,0)$, $\dfrac{\partial u}{\partial x}(t,1)$ exist for all $t \in (0,T]$ and equation (2.1) holds. Suppose that (2.2) and (4.11) hold. Then the following estimate holds for all $\zeta \in [0,\sigma)$, $g_0, g_1 : (0,T] \to (0,+\infty)$, $k_0, k_1 : (0,T] \to [1,+\infty)$ and $t \in (0,T]$:*

$$\|u[t]\|_\infty \leq \max\left(\exp(-\zeta t)\|u[0]\|_\infty, \sup_{0<s\leq t}\left(\max\left(r_0(s), r_1(s), \frac{\|f[s]\|_\infty}{\sigma - \zeta}\right)\exp(-\zeta(t-s))\right)\right) \quad (4.22)$$

*where*

$$r_0(t) := \min\left(|u(t,0)|, \left|g_0(t)\frac{\partial u}{\partial x}(t,0) - k_0(t)u(t,0)\right|\right),\ r_1(t) := \min\left(|u(t,1)|, \left|g_1(t)\frac{\partial u}{\partial x}(t,1) + k_1(t)u(t,1)\right|\right),$$
$$\text{for } t \in (0,T] \quad (4.23)$$

**Proof:** Let (arbitrary) $\zeta \in (0,\sigma)$, $t_1, t_2 \in (0,T)$ with $t_1 < t_2$ be given. Lemma 4.1 implies that the mapping $[t_1, t_2] \ni t \to \|u[t]\|_\infty$ is an absolutely continuous mapping. Moreover, there exists a Lebesgue measure zero set $N \subset [t_1, t_2]$ such that equation (4.13) holds. Moreover, we assume that the Lebesgue measure zero set $N \subset [t_1, t_2]$ has been selected in such a way that (2.1) holds for all $t \in [t_1, t_2] \setminus N$ and for all $x \in (0,1)$ and $u[t] \in C^2((0,1))$ for all $t \in [t_1, t_2] \setminus N$. Lemma 4.2 implies that (4.14) holds with $I(t) = \{ y \in [0,1] : |u(t,y)| = \|u[t]\|_\infty \}$.

Let (arbitrary) $\varepsilon > 0$ be given. Inequality (4.14) guarantees that the following implication holds for all $t \in [t_1, t_2] \setminus N$:

$$\|u[t]\|_\infty \geq \varepsilon + \max\left( r_0(t), r_1(t), \frac{\|f[t]\|_\infty}{\sigma - \zeta} \right) \Rightarrow \frac{d}{dt}\left( \|u[t]\|_\infty \right) \leq \max_{x \in I(t)}\left( \operatorname{sgn}(u(t,x)) \frac{\partial u}{\partial t}(t,x) \right) \quad (4.24)$$

Pick any $t \in [t_1, t_2] \setminus N$. Working exactly as in the proof of Lemma 4.3 we show that implication (4.16) holds.

We next show that the condition $\|u[t]\|_\infty \geq \varepsilon + \max\left( r_0(t), r_1(t), \frac{\|f[t]\|_\infty}{\sigma - \zeta} \right)$ implies that $I(t) \subseteq (0,1)$.

Suppose that $\|u[t]\|_\infty \geq \varepsilon + \max\left( r_0(t), r_1(t), \frac{\|f[t]\|_\infty}{\sigma - \zeta} \right)$ and $0 \in I(t)$. Definition (4.23) implies that $r_0(t) = \left| g_0(t) \frac{\partial u}{\partial x}(t,0) - k_0(t) u(t,0) \right|$ which combined with the fact $\|u[t]\|_\infty > r_0(t)$ gives

$$k_0(t) u(t,0) - \|u[t]\|_\infty < g_0(t) \frac{\partial u}{\partial x}(t,0) < \|u[t]\|_\infty + k_0(t) u(t,0) \quad (4.25)$$

If $u(t,0) = \|u[t]\|_\infty$ then (4.25) in conjunction with the facts that $k_0(t) \geq 1$, $g_0(t) > 0$ gives $\frac{\partial u}{\partial x}(t,0) > 0$, which contradicts the fact that $u$ has a maximum at $x=0$. If $u(t,0) = -\|u[t]\|_\infty$ then (4.25) in conjunction with the facts that $k_0(t) \geq 1$, $g_0(t) > 0$ gives $\frac{\partial u}{\partial x}(t,0) < 0$, which contradicts the fact that $u$ has a minimum at $x=0$. We conclude that the condition $\|u[t]\|_\infty \geq \varepsilon + \max\left( r_0(t), r_1(t), \frac{\|f[t]\|_\infty}{\sigma - \zeta} \right)$ implies that $0 \notin I(t)$.

Suppose that $\|u[t]\|_\infty \geq \varepsilon + \max\left( r_0(t), r_1(t), \frac{\|f[t]\|_\infty}{\sigma - \zeta} \right)$ and $1 \in I(t)$. Definition (4.23) implies that $r_1(t) = \left| g_1(t) \frac{\partial u}{\partial x}(t,1) + k_1(t) u(t,1) \right|$ which combined with the fact $\|u[t]\|_\infty > r_1(t)$ gives

$$-k_1(t) u(t,1) - \|u[t]\|_\infty < g_1(t) \frac{\partial u}{\partial x}(t,1) < \|u[t]\|_\infty - k_1(t) u(t,1) \quad (4.26)$$

If $u(t,1) = \|u[t]\|_\infty$ then (4.26) in conjunction with the facts that $k_1(t) \geq 1$, $g_1(t) > 0$ gives $\frac{\partial u}{\partial x}(t,1) < 0$, which contradicts the fact that $u$ has a maximum at $x=1$. If $u(t,1) = -\|u[t]\|_\infty$ then (4.26) in conjunction with the facts that $k_1(t) \geq 1$, $g_1(t) > 0$ gives $\frac{\partial u}{\partial x}(t,1) > 0$, which contradicts the fact that $u$ has a minimum at $x=1$. We conclude that the condition $\|u[t]\|_\infty \geq \varepsilon + \max\left(r_0(t), r_1(t), \frac{\|f[t]\|_\infty}{\sigma - \zeta}\right)$ implies that $1 \notin I(t)$.

Since the inequality $\|u[t]\|_\infty \geq \varepsilon + \max\left(r_0(t), r_1(t), \frac{\|f[t]\|_\infty}{\sigma - \zeta}\right)$ implies that $I(t) \subseteq (0,1)$, we obtain from (4.11), (4.24) and (4.16) for all $t \in [t_1, t_2] \setminus N$:

$$\|u[t]\|_\infty \geq \varepsilon + \max\left(r_0(t), r_1(t), \frac{\|f[t]\|_\infty}{\sigma - \zeta}\right) \Rightarrow \frac{d}{dt}\left(\|u[t]\|_\infty\right) \leq -\zeta \|u[t]\|_\infty \qquad (4.27)$$

Implication (4.27) and Lemma 2.14 on page 82 in [10] imply that the following estimate holds for all $t \in [t_1, t_2]$:

$$\|u[t]\|_\infty \leq \max\left(\exp(-\zeta(t-t_1))\|u[t_1]\|_\infty, \sup_{t_1 \leq s \leq t}\left(\max\left(r_0(s), r_1(s), \frac{\|f[s]\|_\infty}{\sigma - \zeta}\right)\exp(-\zeta(t-s))\right)\right) + \varepsilon \qquad (4.28)$$

Since $\varepsilon > 0$ is arbitrary, we conclude from (4.28) that the following estimate holds for all $t \in [t_1, t_2]$:

$$\|u[t]\|_\infty \leq \max\left(\exp(-\zeta(t-t_1))\|u[t_1]\|_\infty, \sup_{t_1 \leq s \leq t}\left(\max\left(r_0(s), r_1(s), \frac{\|f[s]\|_\infty}{\sigma - \zeta}\right)\exp(-\zeta(t-s))\right)\right) \qquad (4.29)$$

Finally, since $t_1, t_2 \in (0, T)$ are arbitrary and since $u \in C^0([0,T] \times [0,1])$ (which implies that $\lim_{t_1 \to 0^+}\left(\|u[t_1]\|_\infty\right) = \|u[0]\|_\infty$ and $\lim_{t_2 \to T^-}\left(\|u[t_2]\|_\infty\right) = \|u[T]\|_\infty$), we conclude from (4.29) that (4.22) holds for all $t \in (0, T]$ and $\zeta \in (0, \sigma)$.

Continuity arguments guarantee the fact that estimate (4.22) holds for $\zeta = 0$ as well. The proof is complete. ◁

We are now ready to provide the proof of Theorem 2.2.

**Proof of Theorem 2.2:** Theorem 2.2 is a direct consequence of Lemma 4.4 applied to the function $w \in C^0([0,T] \times [0,1])$ defined by (4.20). ◁

Finally, we provide the proof of Corollary 2.3.

**Proof of Corollary 2.3:** We simply apply Theorem 2.2 with

- $g_0(t) \equiv \frac{\mu_0 \eta(0)}{|\mu_0 \eta'(0) - \lambda_0 \eta(0)|}$, $k_0(t) \equiv 1$ and arbitrary $g_1 : (0,T] \to (0,+\infty)$, $k_1 : (0,T] \to [1,+\infty)$ when (2.9) holds,

- $g_1(t) \equiv \frac{\mu_1 \eta(1)}{\mu_1 \eta'(1) + \lambda_1 \eta(1)}$, $k_1(t) \equiv 1$ and arbitrary $g_0 : (0,T] \to (0,+\infty)$, $k_0 : (0,T] \to [1,+\infty)$ when (2.10) holds,

- $g_0(t) \equiv \frac{\mu_0 \eta(0)}{|\mu_0 \eta'(0) - \lambda_0 \eta(0)|}$, $k_0(t) \equiv 1$, $g_1(t) \equiv \frac{\mu_1 \eta(1)}{\mu_1 \eta'(1) + \lambda_1 \eta(1)}$, $k_1(t) \equiv 1$, when (2.11) holds.

The proof is complete. ◁

# 5. Concluding Remarks

The present paper provided tools for the derivation of ISS estimates in the sup norm of the state. More specifically, we provided novel ISS-style maximum principle estimates which are valid for classical solutions of parabolic PDEs. The obtained results can be applied to highly nonlinear 1-D parabolic PDEs, as illustrated by three illustrative examples.

However, the main results have some disadvantages. The main disadvantages of the main results of the paper are:

1) the fact that they can only be applied to PDEs for which classical solutions can be proved to exist (at least locally; see [3,17]). In some cases, we are in a position to prove the existence of a continuous solution (so that the sup norm estimates make sense) which satisfies the PDE in a weak sense (see for example [25]). For such cases, Theorem 2.1 and Theorem 2.2 cannot be applied. Therefore, there is a need to extend the results for weaker notions of solutions. This is going to be the topic of future research.

2) the fact that they can only be applied to parabolic PDEs with one spatial dimension. The extension to parabolic PDEs with $n-$ dimensional domains is going to be another topic of future research.

**Acknowledgements:** The authors would like to thank Professor Renkun Chen for discussions about the dependence of heat conductivity of solids on temperature.